\documentclass{amsart}
\usepackage{amsthm}
\usepackage{amssymb}
\usepackage{amsmath}
\usepackage{latexsym} 
\newcommand{\sbq}{\subseteq}
\newcommand{\spq}{\supseteq}
\newcommand{\mc}{\mathcal}

\newcommand{\tbf}{\textbf}
\newcommand{\mbf}{\mathbf}
\newcommand{\mf}{\mathfrak}

\newcommand{\inv}{^{-1}}

\newcommand{\bz}{\mbf{0}}
\newcommand{\rad}[1]{\mbf{J}(#1)}
\newcommand{\mat}[2]{\mbf{M}_{#1}(#2)}
\newcommand{\tri}[2]{\mbf{T}_{#1}(#2)}
\newcommand{\B}[1]{\mbf{B}(#1)}
\newcommand{\Sh}[1]{\mc{P}(#1)}
\DeclareMathOperator{\spec}{Spec}

\DeclareMathOperator{\rg}{rank}

\newcommand{\nil}[1]{\mbf{Nil}(#1)}

\newcommand{\U}{\mbf{U}}
\newcommand{\End}[1]{\text{End}(#1)}
\newtheorem{thm}{Theorem}
\newtheorem{lem}[thm]{Lemma}

\newtheorem{prop}[thm]{Proposition}
\newtheorem{obs}[thm]{Observation}
\newtheorem{ex}[thm]{Example}

\newtheorem{cor}[thm]{Corollary}

\newtheorem{defn}[thm]{Definition}
\newtheorem{?}[thm]{Question}

\def\N{\mathbf N}

\def\Z{\mathbf Z}

\newcommand{\varep}{\varepsilon}
\newcommand{\res}{\raisebox{-.5ex}{$|$}}
\begin{document}
\title[Strongly clean matrices]{On strongly clean matrices over commutative clean rings}

\author{Walter Burgess}
\address{University of Ottawa, Ottawa, Canada}
\email{wburgess@uottawa.ca}
\date{\today}
\begin{abstract} The literature about strongly clean matrices over commutative rings is quite extensive. The sharpest results are about matrices over commutative local rings, for example those by Borooah, Diesl and Dorsey. The purpose of this note is to show that, using Pierce sheaf techniques, many of the known results about matrices over commutative local rings can be extended to those over commutative clean rings in general. \end{abstract}
\maketitle

\noindent \tbf{1. Introduction.} Recall that an element $r$ in a ring $R$ is called \emph{clean} if there is an idempotent $e$ and a unit $u$ such that $r=e+u$. If for a clean element $r$,  $e$ and $u$ can be found such that  $eu=ue$ then $r$ is called \emph{strongly clean}. A ring $R$ is called \emph{(strongly) clean} if each of its elements is (strongly) clean. Clean rings are special cases of exchange rings and the two classes coincide for rings whose idempotents are all central (\cite{N77} and \cite{Cetal06}). 

Strongly clean matrices over commutative rings have been extensively studied. Some of the sharpest results are known about matrices over commutative local rings (\cite{Betal08} and its list of references).  The purpose of this note is to show, using Pierce sheaf techniques, how, in many instances, the results from the commutative local case can be extended to matrices over all commutative clean rings.  Some recent extensions of the results about individual strongly clean matrices over commutative local rings to commutative rings such that finitely generated projective modules are free are also viewed from a Pierce point of view. However, an example shows some of the limits to further generalization even to Dedekind domains.\\

\noindent\tbf{2. Notation, terminology and background.} Throughout rings are unital.  If $R$ is a ring then $\mat{n}{R}$ denotes the ring of $n\times n$ matrices over $R$; $\U(R)$ is the group of units of $R$ and $\rad{R}$ is the Jacobson radical. For $R$ commutative and $A\in \mat{n}{R}$, $\chi(A)$ denotes the characteristic polynomial of $A$.

The following lemma is essential and here will be applied to $A\in \mat{n}{R}$ but is true in all endomorphism rings of modules. 

\begin{lem}\cite[Theorem~3]{N99}, \cite[Proposition~2.3]{Cetal06} \label{dsdecomp} For any ring $R$ and $R$-module $M_R$, $\varphi\in \End{M}$ is strongly clean if and only if $M$ has a $\varphi$-invariant decomposition $M= A\oplus B$ where $\varphi\res_A\colon A\to A$ is an automorphism and $(1-\varphi)\res_B\colon B\to B$ is an automorphism. \end{lem}

The following definition comes from \cite[Definition~5]{Betal08}.
\begin{defn} \label{SR} Let $R$ be a commutative ring and $f\in R[t]$ a monic polynomial. A factorization $f=f_0f_1$ is called an SR-\emph{factorization} if $f_0(0)$ and $f_1(1)$ are in $\U(R)$. If, in addition, the ideal $(f_0,f_1) = R[t]$, the factorization is called an SRC-\emph{factorization}.\end{defn}

This form of the definition is equivalent to that given in \cite{Betal08} as noted in that article. The choice of terminology is explained in \cite[page~283]{Betal08}.

This note intends to exploit the following theorem from \cite{Betal08} which, in turn, uses basic results from \cite{N99}.
\begin{thm} \cite[Theorem~12]{Betal08} \label{Betal1}Suppose $R$ is a commutative local ring and $n\ge 2$, $h\in R[t]$ is a monic polynomial of degree $n$. Then the following are equivalent:

(1) For any $A\in \mat{n}{R}$ with the characteristic polynomial $\chi(A)=h$, $A$ is strongly clean.

(2) The companion matrix $C_h$ of $h$ is strongly clean.

(3) There exists an SRC-factorization of $h$ in $R[t]$.
\end{thm}

The methods of the Pierce sheaf are ideal in this situation. See \cite[Chapter~5 (2)]{J82} or \cite{BS76} for details, both based on \cite{P67}. The reason for this is that all the statements in Theorem~\ref{Betal1} can be expressed in terms of elements of the ring $R$ using a finite number of existential quantifiers and a finite number of equations. The next section sketches, in very summary form, the key elements of Pierce sheaf theory which will be needed. Everything mentioned can be done without using sheaf terminology but it is intuitively natural to use it especially in this context.\\

\noindent
\tbf{3. Pierce sheaf background material.} Throughout, $R$ will stand for a commutative ring and $\B{R}$ will be its set of idempotents, viewed either as a boolean ring or as a boolean algebra (according to the addition chosen).  The spectrum of $\B{R}$, $\spec \B{R}$, is a compact totally disconnected Hausdorff space.  The basic closed sets of $\spec \B{R}$ are of the form, for $e\in \B{R}$, $D(e)=\{x\in \spec \B{R}\mid e\notin x\}$. Note that these sets are both closed and open.  A \emph{complete orthogonal set of idempotents} is a (finite) orthogonal set of idempotents from $\B{R}$ whose sum is 1.

The \emph{Pierce sheaf} for $R$, denoted $\Sh{R}$, has $\spec \B{R}$ as base space and the stalks are the rings $R_x=R/xR$, for $x\in \spec \B{R}$. For $r\in R$ and $x\in \spec \B{R}$, $r+xR$ is written $r_x$ and if $V\sbq R$ is a subset then the set of images in $R_x$ is $V_x$. Note that these stalks have no non-trivial idempotents. The disjoint union of these stalks is the \emph{espace \'etal\'e} and has basic open sets of the form, for $r\in R$ and $D(e)\sbq \spec \B{R}$, $\{r_x\mid x\in D(e)\}$. There are two key statements, the second follows from the first by exploiting the compactness of $\spec \B{R}$ and the fact that a finite set of idempotents $\{e_1, \ldots, e_k\}$ such that $\bigcup_{1=1}^kD(e_i) = \spec \B{R}$ can be transformed into a complete orthogonal set by subtractions and multiplications.
\begin{thm} \cite[Theorem page 183 and its proof]{J82} \label{basic} Let $R$ be a commutative ring. Then:

(1) if for $r, s\in R$ and $x\in \spec \B{R}$, $r_x=s_x$ then there is $e\in \B{R}$ such that $x\in D(e)$ and for all $y\in D(e)$, $r_y=s_y$. In other words, $re=se$.

(2) the ring of global sections of $\Sh{R}$ is isomorphic to $R$.
\end{thm}

In fact the ring $R$ need not be commutative but, if not, $\B{R}$ needs to be taken as the set of \emph{central} idempotents. 

The key method used below, which is also the tool used to prove that (1)~$\Rightarrow$~(2) in Theorem~\ref{basic}, is the following.  Suppose for $a,b\in R$ that for each $x\in \spec \B{R}$ the equation $a_x=b_x$ is true. Then, by Theorem~\ref{basic}~(1), for each $x\in \spec \B{R}$ there is $e(x)\in \B{R}$ with $e(x)\notin x$ so that $ae(x)=be(x)$.  Then, $\bigcup_{x\in\spec \B{R}} D(e(x)) =\spec \B{R}$. Compactness reduces this cover to a finite cover, say $\{D(e_1), \ldots, D(e_k)\}$; the set of idempotents $\{e_1,\ldots, e_k\}$ can be transformed  into a complete orthogonal set of idempotents $\{\varep_1,\ldots, \varep_n\}$. Then, $a = \sum_{i=1}^n a\varep_i =\sum_{i=1}^n b\varep_i =b$. Notice that each $D(\varep_i)\sbq D(e_j)$, for some $e_j$. This process lifts an equation which holds at each of the stalks to an equation in $R$. This extends readily to any finite set of such equations. 

One more fact is required (which does not require commutativity).

\begin{lem}\cite[Theorem~3.4]{BS79}, \cite[Proposition~1.8(2)]{N77} and \cite[Proposition page~187]{J82} \label{localstalks} The following are equivalent for a ring $R$.

(1) All the Pierce stalks of $R$ are local rings.

(2) $R$ is an exchange ring all of whose idempotents are central. 

(3) $R$ is a clean ring all of whose idempotents are central.\end{lem}

 In particular, Lemma~\ref{localstalks} characterizes commutative clean rings as those whose Pierce stalks are local.\\
 
\noindent\tbf{4. The results.}
The first step is to show how Theorem~\ref{Betal1} can be generalized to all commutative clean rings. However, the concept of an SR(C)-factorization needs to be adjusted.

\begin{defn} \label{src-gen} Let $R$ be a commutative ring and $h\in R[t]$ a monic polynomial of degree $n\ge 1$. Then, $h$ is said to have a  gSR-\emph{factorization} (gSRC-\emph{factorization}) if there is a complete orthogonal set of idempotents in $R$, depending on $h$, $\{e_1,\ldots ,e_k\}$ such that for $i=1,\ldots, k$, $he_i$ has an SR(C)-factorization in $e_iR[t]$. \end{defn}

The ``g'' may be read as ``globalized'' or ``generalized'' according to taste. The following simple example illustrates the difference between SR-factorizations and gSR-factorizations. Consider $R= \Z_{(2)}\times \Z_{(2)}$ and the polynomial $h=t^2+(3,1)t+ (2,3)$; it has no SR-factorization but does have a gSRC-factorization.  Theorem~\ref{main}, below, or two applications of Theorem~\ref{Betal1} will show that any $A\in \mat{n}{R}$ with $\chi(A)= h$ is strongly clean.

\begin{thm}\label{main} Suppose $R$ is a commutative clean ring, $n\ge 2$ and $h\in R[t]$ is a monic polynomial of degree $n$. Then, the following are equivalent: 

(1) For any $A\in \mat{n}{R}$ with the characteristic polynomial $\chi(A)=h$, $A$ is strongly clean.

(2) The companion matrix $C_h$ of $h$ is strongly clean.

(3) There exists a gSRC-factorization of $h$ in $R[t]$.
  
\end{thm}

Proof. For $A\in \mat{n}{R}$ and $x\in V$, let the matrix in $\mat{n}{R_x}$ obtained by reducing the entries modulo $xR$ be denoted $A_x$. Similarly, for $h\in R[t]$, its reduction modulo $xR$ is called $h_x$. Clearly $C_{h_x}=(C_h)_x$. Moreover, an SRC-factorization of $h$ in $R[t]$ reduces modulo $xR$ to an SRC-factorization of $h_x$ in $R_x[t]$. (Keep in mind that $(h_0,h_1)=R[t]$ can be expressed by the existence of $u,v\in R[t]$ with $h_0u+h_1v=1$.)

If $A\in \mat{n}{R}$ is strongly clean, then, by reduction modulo $xR$, so is $A_x$. 

The conclusion is that if (1), (2) or (3) holds over $R$ then, for each $x\in \spec\B{R}$ the corresponding condition holds over $R_x$. However, each $R_x$ is a local ring and there all three conditions are equivalent by Theorem~\ref{Betal1}. 

Now suppose that statements (1), (2) and (3) are true over $R_x$ for all $x\in \spec\B{R}$ for a given matrix $A$ and $h\in R[t]$ of degree $n$. The key is to note that each of the conditions in $R_x$ can be expressed by statements involving finitely many existential quantifiers and finitely many equations over $R_x$. For example (3) becomes $\exists g_{0,x}, g_{1,x}\in R_x[t]$ such that $g_{0,x}(0_x)=h_x$ and $g_{1,x}(1_x)$ are in $\U(R_x)$ and $\exists u(x),v(x)\in R_x[t]$ such that $g_{0,x}u(x)+ g_{1,x}v(x)=1_x$. For a fixed $x\in V$ these statements become a finite set of equations over $R_x$ using just the coefficients. However, this implies that there is $e(x)\in \B{R}$, $e(x)\notin x$ such that when all the coefficients in all the equations are lifted to elements of $R$, say to form a set of elements $T\sbq R$; multiplying them by $e(x)$ yields the same equations holding modulo $y$ for all $y\in D(e(x))$.  Then the set of elements $e(x)T\sbq e(x)R$ satisfy all the equations dealt with here modulo $yR$ for each $y\in D(e(x))$.  In other words,  $e(x)h$ has an SRC-factorization in the ring $e(x)R[t]$. 

The next step is convert the cover of $\spec\B{R}$ given by the $D(e(x))$ to one given by a complete finite orthogonal set of idempotents, say $\{\varep_1,\ldots, \varep_m\}$. Each of the sets $D(\varep_i)$  is inside one of the $D(e(x))$ and it follows that $\varep_ih$ has an SRC-factorization in $\varep_iR[t]$.  Label all the resulting polynomials as follows: $\varep_ih= \varep_ih_{0,i}h_{1,i}$ where $h_{0,i}(0), h_{1,i}(\varep_i)\in \U(\varep_iR)$. Moreover, there are $u_i,v_i\in \varep_iR[t]$ such that $h_{0,i}u_i + h_{1,i}v_i=\varep_i$.  This gives statement~3.

The other two statements lift in the same manner. For example if (1) is true for $A_x$ for all $x\in \spec\B{R}$ then the following are true at $x$: there exist $U(x), V(x),E(x)\in \mat{n}{R_x}$ so that $U(x)V(x)=V(x)U(x)= I_n$ (the identity matrix in $\mat{n}{R_x}$), $E(x)^2=E(x)$, $A_x= U(x)+E(x)$ and $U(x)E(x)= E(x)U(x)$.  When all this is expressed in terms of the coefficients there is again a finite set of equations over $R_x$.  The lifting process proceeds in the same manner. The finite orthogonal set of idempotents need not be the same as before, although both parts could have been done simultaneously, but since the resulting equations are over $R$ this does not matter. 

Of course, statement (2) lifts as well. 

Hence, if one of the three statements holds over $R_x$ for each $x\in \spec\B{R}$ then all three statements hold over $R$.  It follows that if any one of the three conditions hold over $R$ then all three do. \qed \\[-.5ex]

\begin{cor}\label{deg} Suppose in Theorem~\ref{main} the three conditions are satisfied for $h\in R[t]$ of degree $n$. Then the complete orthogonal set of idempotents $\{\varep_1,\ldots ,\varep_m\}$ showing $h$ has a gSRC-factorization may be chosen with $m\le n+1$.  \end{cor}

Proof. In the proof of Theorem~\ref{main} once the $\varep_i$ have been found it can be seen that the degree of $(h_{0,i})_x$ is the same for all $x\in D(\varep_i)$.  In fact once the SRC-factorizations of $h_x= g_{0,x}g_{1,x}$ have been chosen, it is easy to see that the function $\tau \colon \spec \B{R}\to \N\cup \{0\}$ with $\tau(x) = \deg g_{0,x}$, is continuous when $\N\cup \{0\}$ has the discrete topology.  Define $\zeta_j\in \B{R}$ by $D(\zeta_j)= \tau\inv(j)$; then, ignoring those which are 0, $\{\zeta_0,\ldots, \zeta_n\}$ is a complete orthogonal set of idempotents.  For each $j$, $0\le j\le n$, it follows that $\zeta_j$ is a sum of some of the $\varep_i$. It then suffices to sum the components making up the statement that $\varep_ih=h_{0,i}h_{1,i}$ is an SRC-factorization in $\varep_iR[t]$ over all $i$ where $\tau(x)=j$ for $x\in D(\varep_i)$. \qed \\[-.5ex]

In \cite[Corollary~15]{Betal08}, the authors give five equivalent conditions for the \emph{ring} $\mat{n}{R}$ to be strongly clean, where $R$ is a commutative local ring. All these conditions can be globalized but only two of them will be noted here.

\begin{cor}\label{allstrcl} Let $R$ be a commutative ring and $n\in \N$, $n\ge 2$. Then the following are equivalent. 

(1) $\mat{n}{R}$ is a strongly clean ring.

(2) For every monic $h\in R[t]$ of degree $n$, $h$ has an gSRC-factorization.  \end{cor}

Proof. Since ``strongly clean'' is preserved in corner rings (\cite[Theorem~2.4]{C06}), if (1) holds then $R$ is (strongly) clean.  Hence Theorem~\ref{main} applies showing that (2) holds. 

Assume (2).  Fix $x\in \spec\B{R}$; it is first shown that $R_x$ is local.  Suppose $b\in R$ with $b_x\notin \U(R_x)$.  Consider $h=(t^2+b)t^{n-2}$. Since $h$ has a gSR-factorization, say using the complete orthogonal set of idempotents $\{e_1,\ldots ,e_k\}$, suppose $x\in D(e_i)$ and $he_i=f_{0,i}e_if_{1,i}e_i$ is a gSRC-factorization; it follows that $(f_{0,i})_x$ is a factor of $t^2+b_x$.  If $(f_{0,i})_x$ is the trivial factor $1_x$ then $(f_{1,i})_x=h_x$ and $(f_{1,i})_x(1_x) =1_x+b_x\in \U(R_x)$. Otherwise, $t^2+b_x = (t+c_x)(t+d_x)$ where $c_x\in \U(R_x)$ and $c_xd_x=b_x$. However, $c_x+d_x=0$ which would show that $d_x$ and $c_xd_x=b_x$ were units. This is not possible. Hence, $1_x+b_x\in \U(R_x)$.  Since this applies to $b_x$ it also applies to $b_xr_x$ for all $r\in R$. Since every non-unit of $R_x$ is in $\rad{R_x}$ it follows that  $R_x$ is a local ring. 

Now Theorem~\ref{basic} applies (since these are, by assumption, SRC-factor\-izations) showing that $\mat{n}{R_x}$ is strongly clean. Pierce techniques then make $\mat{n}{R}$ strongly clean.
\qed\\[-.5ex]


The Pierce sheaf theory also applies to modules over a (commutative) ring $R$.  If $M$ is an $R$-module then the $R_x$-modules $M_x= M/xM$ form a sheaf of modules over $\spec \B{R}$ in much the same way as constructing the sheaf for the ring itself. There is a result analogous to Theorem~\ref{basic}.  

In \cite{CKK13} the authors generalize some of the results of \cite{Betal08} to commutative rings $R$ so that every finitely generated projective module is free, called \emph{projective free rings}. Such rings have no non-trivial idempotents but do not need to be local; indeed, all PIDs are examples. Hence, Theorem~\ref{main} has an immediate corollary.
\begin{cor}\label{main1} Let $R$ be a commutative ring such that its Pierce stalks are all projective free rings. Then the three statements of Theorem~\ref{main} are equivalent for $R$.  \end{cor}

Proof. Theorem~2.4 of \cite{CKK13} gives the equivalence of the three statements for projective free rings and the proof of Theorem~\ref{main} shows how this lifts to the ring $R$. \qed \\[-.5ex]

It would be convenient to have more information about  the rings of Corollary~\ref{main1}. A preliminary lemma, which probably exists in the literature, is useful.

\begin{lem} \label{liftpr} Let $R$ be a commutative ring and for some $x\in \spec \B{R}$ let $P$ be a finitely generated projective $R_x$-module. Then, there exists a finitely generated projective $R$-module $V$ such that $V_x=P$.  \end{lem}

Proof. Suppose that $P$ is a summand of $(R_x)^m$. In $R^m$, fix the standard basis $\{\kappa_1, \ldots, \kappa_m\}$; then $\{(\kappa_1)_x, \ldots, (\kappa_m)_x\}$ is the standard basis for $(R_x)^m$ which is used to give projective coordinates $$\{p_1, \ldots , p_m; \lambda_1, \ldots, \lambda_m\}$$ for $P$.  Now lift the elements $p_i$ to $b_i\in R^m$, i.e., $(b_i)_x=p_i$, $i=1, \ldots, m$. Put $M= \sum_{i=1}^mb_iR$. Denote the projections $R^m\to R$ using the standard basis by $\tau_i$, $i=1, \ldots, m$.  The key equations are $$ \sum_{i=1}^m p_i\lambda_i(p_j) = p_j,\;\;  j=1, \ldots, m\;. $$
These translate to $$ (\sum_{i=1}^m b_i\tau_i(b_j))_x = (b_j)_x, \;\; j= 1,\ldots ,m \;.  $$ Hence, there is $e\in \B{R}\setminus x$ with $(\sum_{i=1}^m b_i\tau_i(b_j))e = b_je, j=1, \ldots ,m$.  It follows that $\{b_1e, \ldots, b_me; \tau_1e,\ldots, \tau_me\}$ is a system of projective coordinates for $Me$ as an $Re$-module. Now set $V= R^m(1-e) \oplus Me$. Note that $V_x=P$. 

Let $\varphi \colon (Re)^m\to Me$ be given by $\varphi (\varep_ie) = b_ie$, $i=1,\ldots, m$. This splits via, say, $\psi \colon Me\to (Re)^m$. Define $\zeta \colon R^m\to V$ by $\zeta (b(1-e) + me)= b(1-e)+\varphi (me)$, $b\in R^m$ and $m\in M$. This is split by the identity on $R^m(1-e)$ and $\psi $ on $Me$. Hence, $V$ is a finitely generated projective $R$-module. \qed

\begin{prop}\label{prfree}  Let $R$ be a commutative ring.  All the Pierce stalks of $R$ are projective free rings if and only if for every finitely generated projective $R$-module $P$ there is a complete orthogonal set of idempotents in $R$, $\{e_1,\ldots, e_k\}$, the idempotents and $k$ depending on $P$, such that for $1\le i\le k$, $Pe_i$ is a free $Re_i$-module. \end{prop}

Proof.  The discussion in \cite[\S 13]{P67} is about modules over commutative regular rings but readily adapts to the present situation. A (split) short exact sequence of modules $\bz \to A\to B\to C\to \bz$ gives rise, for each $x\in \spec \B{R}$, to a (split) short exact sequence $\bz \to A_x\to B_x \to C_x \to\bz$. 

Now suppose that all the stalks are projective free.  Let $P$ be a finitely generated projective $R$-module with $P\oplus Q\cong R^n$, some $n$. Let $\{a_1,\ldots, a_n\}$ be a set of generators for $P$. For each $x\in \spec \B{R}$, $P_x$ is free of rank, say $k(x)$. Since each $(a_i)_x$ can be expressed in terms of $k(x)$ generators in $P_x$; using these expressions there is an idempotent $e(x)\notin x$ so that each $a_ie(x)$ can be expressed in terms of $k(x)$ generators.  Hence, $U=\{y\in \spec \B{R}\mid \rg P_y\le k(x)\}$  is an open neighbourhood in $\spec \B{R}$. 

The next step is to consider the splitting $P\oplus Q\cong R^n$. For $x$ as in the previous paragraph, $\rg (P_x) + \rg (Q_x) = n$. It follows that $V=\{y\in \spec\B{R} \mid \rg (Q_y)\le n-k(x)\}$ is also open. However, $V$ is the complement of $U$, showing that $U$ is also closed. From this it can be concluded that the function $\rho\colon \spec \B{R} \to \{0, 1, \ldots, n\}$ given by $\rho (x) = \rg (P_x)$  is continuous. Put $e_i\in \B{R}$ to be the idempotent with $\rho\inv (i)= D(e_i)$.  The set $\{e_0,\ldots, e_n\}$ (ignoring zeros) is a complete orthogonal set of idempotents and for each $i$, $Pe_i$ is a free $Re_i$-module, as required. 

For the converse, fix $x\in \spec \B{R}$ and suppose that $P$ is a finitely generated projective $R_x$-module. According to Lemma~\ref{liftpr}, there is a finitely generated projective $R$-module $V$ such that $V_x=P$. By hypothesis, there is a complete orthogonal set of idempotents $\{e_1,\ldots ,e_k\}$ such that $Ve_i$ is a free $Re_i$-module, $i=1,\ldots, k$. There is a unique $j$ such that $e_j\notin x$. Then $P= (Ve_j)_x$ is a free $R_x$-module.
\qed \\[-.5ex]
 
Recall that an element $r\in R$ is \emph{strongly $\pi$-regular} if both chains of one-sided ideals  $rR\spq r^2R \spq r^3R \spq  \cdots$ and $ Rr\spq Rr^2 \spq 
Rr^3 \spq \cdots $ terminate.  If $A\in \End{M_R}$ then $A$ is strongly $\pi$-regular if and only if there is a decomposition $M=P\oplus Q$ into $A$-invariant submodules where $A\res_P$ is an automorphism and $A\res_Q$ is a nilpotent endomorphism; hence, in particular, a strongly $\pi$-regular element is strongly clean.  (See \cite[page~3589]{N99} and \cite[Lemma~40 and footnote page~292]{Betal08}.)

There is a characterization \cite[Proposition~44]{Betal08} of when a matrix over a commutative local ring is strongly $\pi$-regular. It involves another sort of factorization of the characteristic polynomial.  The symbol $\nil{R}$ is the nil radical of a commutative ring $R$. Notice that for any commutative ring $R$,  $a\in \nil{R}$,  if and only if for all $x\in \spec \B{R}$, $a_x\in \nil{R_x}$.  

\begin{defn} \label{SP} \cite[Definitions~42 and 43]{Betal08}  Let $h\in R[t]$ be a monic polynomial, $R$ a commutative ring. A factorization $h=h_0p_0$ is called an SP-\emph{factorization} if $h_0(0)\in \U(R)$ and $p_0-t^{\deg(p_0)} \in \nil{R}[t]$. Just as in Definition~\ref{src-gen}, this is generalized (globalized) to a gSP-\emph{factorization}.\end{defn}


Using the decomposition of \cite[Lemma~40]{Betal08}, it is shown in \cite[Proposition~44]{Betal08} that the following statements are equivalent for a commutative local ring $R$ and a monic $h\in R[t]$: (1)~Every $A\in \mat{n}{R}$ with $\chi(A)=h$ is strongly $\pi$-regular, (2)~There exists $A\in \mat{n}{R}$ with $\chi(A)=h$ is strongly $\pi$-regular, and (3)~$h$ has an SP-factorization. 

\begin{prop} \label{pi-reg} Let $R$ be a commutative clean ring. Then the following are equivalent for a monic $h\in R[t]$ of degree $n$: 

(1)~Every $A\in \mat{n}{R}$ with $\chi(A)=h$ is strongly $\pi$-regular.

(2)~There exists $A\in \mat{n}{R}$ with $\chi(A)=h$ which is strongly $\pi$-regular.

(3)~The polynomial $h$ has a gSP-factorization. \end{prop} 

Proof.  If, in a ring $S$, some $a\in S$ is strongly $\pi$-regular with, say $a^{k+1}= s^kx$ and $a^{l+1} = ya^l$ then for any $m\ge \max(k,l)$ $a^{m+1} = a^mx = ya^m$. 
This simple remark means that if $a\in S$ is such that, for all $x\in \spec \B{S}$, $a_x$ is strongly $\pi$-regular, then by the methods already used it is possible to conclude that $a$ is strongly $\pi$-regular. 

In the proposition, (2) is a special case of (1).  If, for some $A\in \mat{n}{R}$ with $\chi(A)=h$, $A$ is strongly $\pi$-regular there is a decomposition $R^n=P\oplus Q$ with $A\res_P$ an isomorphism and $A\res_Q$ is nilpotent. According to Proposition~\ref{prfree} applied twice to $P$ and to $Q$, there is a complete orthogonal set of idempotents $\{e_1,\ldots, e_k\}$ such that $Pe_i$ and $Qe_i$ are free $Re_i$-modules, $i=1,\ldots ,k$. Notice that for all $x\in D(e_i)$,  $\rg P_x$ is constant as is $\rg Q_x$.  Then for each $x\in D(e_i)$, $h_x$ has an SP-factorization $h_x= h_{0,x}p_{0,x}$ where the degrees of the factors do not depend on $x\in D(e_i)$. This is because the proof of \cite[Proposition~41]{Betal08} shows that $\deg p_{0,x}= \rg Q_x$. Then, refining the complete orthogonal set of idempotents if necessary, it may be assumed that $he_i$ has an SP-factorization in $Re_i[t]$. This shows (2) $\Rightarrow$ (3).

Finally, assume (3) holds for $h\in R[t]$ and $\chi(A)=h$. For any $i=1,\ldots, k$ and $x\in D(e_i)$, $h_x$ has an SP-factorization and, hence, $A_x$ is strongly $\pi$-regular, showing that $A$ is strongly $\pi$-regular. \qed\\[-.5ex]

As in Corollary~\ref{deg} it can be seen that the number of elements in the complete orthogonal set of idempotents in the gSP-factorization in Proposition~\ref{pi-reg} can be taken to be $\le n+1$ where $n=\deg (h)$.

In \cite[Lemma~2.3]{CKK13} the authors use the resultant of two polynomials to show, over a commutative ring $R$, for two monic polynomials $f, g\in R[t]$ that $(f,g)=R[t]$ if and only if for every maximal ideal $M$ of $R$, $(f,g)$ is not contained in $M[t]$. Proposition~\ref{pi-reg} is proved using Pierce techniques and \cite[Proposition~44]{Betal08}. The key ingredients of \cite[Proposition~44]{Betal08}, which is about local rings, are the fact that a local ring is projective free and that an SP-factorization  (which is clearly an SR-factorization) is, in fact a SRC-factorization. However, the quoted lemma from \cite{CKK13} makes it clear that if $R$ is a commutative projective free ring and $h=h_0p_0$ is an SP-factorization of a monic $h\in R[t]$ then $(h_0,p_0)=R[t]$.  The result yields a more general form of  Proposition~\ref{pi-reg}.

\begin{prop}\label{pi-reg+} Let $R$ be a commutative ring whose Pierce stalks are all projective free rings. Then the following are equivalent for a monic $h\in R[t]$ of degree $n$. 

(1)~Every $A\in \mat{n}{R}$ with $\chi(A)=h$ is strongly $\pi$-regular.

(2)~There exists $A\in \mat{n}{R}$ with $\chi(A)=h$ which is strongly $\pi$-regular.

(3)~The polynomial $h$ has a gSP-factorization. \end{prop}

In \cite{Betal07} the authors show (\cite[Corollary~10]{Betal07}), among many other things, that if $R$ is a commutative local ring then the ring of (say upper) triangular matrix rings over $R$, $\tri{n}{R}$, is strongly clean. The following observation is clear based on the kinds of calculations done frequently above; it does give, in the easy commutative case, a converse to \cite[Lemma~2(2)]{Betal07}.

\begin{obs}\label{triag} If $R$ is a commutative clean ring then for any $n\ge 2$, $\tri{n}{R}$ is strongly clean. \end{obs} 

In \cite[Proposition~30]{Betal08} the authors give a supplementary equivalent condition about when $\mat{2}{R}$ is strongly clean, where $R$ is a commutative local ring. It  is: (3)~The polynomial $t^2-t+a\in R[t]$ has a root (in $R$) for every $a\in \rad{R}$. This type of condition does not readily globalize since, even for a commutative clean ring $R$, it does not follow that $\rad{R_x} = (\rad{R})_x$, for $x\in \spec \B{R}$.  A characterization of rings which this does occur is found in \cite[Theorem~4.1]{BR09}. Three of the conditions given are reproduced here in the commutative case.

\begin{lem} \label{radicals}  Let $R$ be a commutative ring. The following are equivalent.

(1) $R$ is clean and for all $x\in \spec \B{R}$, $\rad{R_x}= (\rad{R})_x$.

(2) For all $r\in R$ there is $e\in \B{R}$ such that $re +(1-e)\in \U(R)$ and $r(1-e)\in \rad{R}$. 

(3) $R/\rad{R}$ is von Neumann regular and idempotents lift uniquely modulo $\rad{R}$. \end{lem} 

Rings satisfying the above three conditions are called J-\emph{clean rings}.  Examples are found at the start of \cite[\S4.2]{BR09}.  Clearly commutative local rings are J-clean.

\begin{prop} \label{jclean} Let $R$ be a commutative J-clean ring. Then the following are equivalent.
 
(1) $\mat{2}{R}$ is strongly clean. 

(2) The polynomial $t^2-t +a\in R[t]$ has a root in $R$ for each $a\in \rad{R}$. \end{prop}

Proof. Assume (1). Given a polynomial $h=t^2-t+a$, $a\in \rad{R}$  it follows from \cite[Proposition~30]{Betal08} that $h_x$ has a root in $R_x$ for each $x\in \spec \B{R}$. If, for a fixed $x\in \spec\B{R}$, $h_x(r_x)=0_x$ then there is $e\in \B{R}\setminus x$ so that $re$ is a root of $he$ in $Re$. The now familiar process gives a root for $h$ in $R$. 

Assume (2). Since $R$ is a J-clean ring, every element of $\rad{R_x}$ is the image of an element of $\rad{R}$. Hence condition (2) holds in $R_x$ for each $x\in \spec \B{R}$. Then, \cite[Proposition~30]{Betal08} shows that each $\mat{2}{R_x}$ is strongly clean. It follows that $\mat{2}{R}$ is strongly clean. \qed  \\[-.5ex]

Two additional equivalent conditions are added in \cite[Proposition~20]{Betal08} in the case where $2\in \U(R)$ and $R$ is local. Assuming that the J-clean ring of Proposition~\ref{jclean} is such that $2\in \U(R)$, it can then be added that conditions (1) and (2) are equivalent to each of (3)~every element of $1+\rad{R}$ has a square root in $1+\rad{R}$ and (4)~every element of $1+\rad{R}$ has a square root in $R$.

The requirement in several of the above results that the commutative ring $R$ be projective free or have stalks projective free does play a key role. In \cite[Corollary~19]{FY08} it is asserted that if $R$ is a commutative projective free ring and $A\in \mat{n}{R}$ is strongly clean then the characteristic polynomial $\chi(A)$ has an SR-factorization (see also \cite[Lemma~9]{Betal08} for the local case). This conclusion can fail when $R$ is a commutative domain which is not projective free.

\begin{ex}\label{notpf} There is a Dedekind domain $R$ which is not projective free and $A\in \mat{2}{R}$ which is strongly clean but whose characteristic polynomial does not have a SR-factorization.\end{ex} Consider the ring $R=\Z[\theta]$ where $\theta=\sqrt{-5}$ and the ideal $\mf{A}=(2,1+\theta)$. This ring and the ideal are examined in detail in \cite[Exercise~2.9]{L07} and the notation established there is used. It is shown that $M=\mf{A}\oplus \mf{A}\cong R\times R$ while $\mf{A}$ is projective but not free. Elements of $M$ are denoted in the form $[a,b]$, $a,b\in \mf{A}$. Then, quoting \cite{L07}, $\{[-2,1-\theta], [1+\theta, -2]\}=\{f_1,f_2\}=\mf{b}$ is a free basis for $M$. Consider $\varphi\in \End{M}$ defined by $\varphi ([a,b]) =[a+b, 2b]$. 

It will first be shown that $\varphi$ is strongly clean using Lemma~\ref{dsdecomp}. Write $M=\{[a,0]\mid a\in \mf{A}\} \oplus \{[b,b]\mid b\in \mf{A}\}= X\oplus Y$. Then, $\varphi\res_X$ is the identity while $(1-\varphi)\res_Y$ is minus the identity. 

The following formulas are used (\cite[page~31]{L07}): 
\begin{alignat*}{4}
[2,0]&=\; &2f_1 + (1-\theta)f_2 &\qquad [1+\theta,0]&=\; &(1+\theta) f_1+3f_2\\
[0,2]&=\; &(1+\theta)f_1 + 2f_2& \qquad [0,1+\theta] &=\; &(-2+\theta)f_1 + (1+\theta)f_2 \notag
 \end{alignat*} 

With this notation, \begin{alignat*}{2}\varphi(f_1) &= [-1-\theta, 2-2\theta]&= (5-\theta) f_1+(-1-2\theta)f_2 \\
\varphi(f_2) &= [-1+\theta, -4] &= (-3-\theta)f_1+(-2+\theta)f_2 \end{alignat*}  

Then the matrix for $\varphi$ with respect to $\mf{b}$ is $$A= \left(\begin{matrix} 5-\theta& -1-2\theta \\ -3-\theta& -2+\theta\end{matrix}\right)$$ and its characteristic polynomial is $\chi(A)= t^2-3t +2-8\theta$ with discriminant $1-32\theta$, whose norm is not a square integer. Thus $\chi(A)$ does not factor in $R[t]$ and neither $\chi(A)(0)$ nor $\chi(A)(1)$ is a unit.  Hence, $\chi(A)$ does not have an SR-factorization.

\end{document}